\documentclass[11pt]{article}
\usepackage[margin=0.85in]{geometry}
\usepackage{type1cm}       
\usepackage{makeidx}
\usepackage{graphicx}        
\graphicspath{{author/figures/}}   
\usepackage{multicol}        
\usepackage[bottom]{footmisc}
\usepackage{csquotes}      
\usepackage{newtxtext}
\usepackage[varvw]{newtxmath}
\usepackage{hyperref}
\usepackage{algorithm,algpseudocode}
\usepackage{pgfplotstable}
\usepackage{array,multirow,graphicx}
\usepackage{pgfplots}
\pgfplotsset{compat=newest}
\usepgfplotslibrary{groupplots}
\usetikzlibrary{matrix}
\DeclareMathAlphabet{\pazocal}{OMS}{zplm}{m}{n}    
\pgfplotsset{compat = 1.3}
\usepackage{wrapfig}
\usepackage{array}
\newcommand{\R}{\mathbb{R}}
\algnewcommand\algorithmiconput{\textbf{Constants:}}
\algnewcommand\algorithmicinput{\textbf{Input:}}
\algnewcommand\algorithmicoutput{\textbf{Output:}}
\usepackage{xcolor}
\definecolor{TR}{HTML}{6a4c93}
\definecolor{SN}{HTML}{c9182c}
\usepackage{mathtools}
\mathtoolsset{showonlyrefs}
\usepackage{booktabs}

\DeclareMathOperator*{\argmax}{arg\,max}

\algnewcommand\Constants{\item[\algorithmiconput]}

\makeindex            
                    
\algtext*{EndFor}
\algtext*{EndWhile}

\begin{document}

\title{A Non-Monotone Preconditioned Trust-Region Method for Neural Network Training}
\author{Andrea Angino\thanks{UniDistance Suisse; andrea.angino@unidistance.ch} \and Bindi Çapriqi\thanks{King Abdullah University of Science and Technology; bindi.capriqi@kaust.edu.sa} \\
Shega Likaj\thanks{Universit\`a della Svizzera italiana; shega.likaj@kaust.edu.sa} \and Ken Trotti\thanks{King Abdullah University of Science and Technology; ken.trotti@kaust.edu.sa} \and Rolf Krause\thanks{Universit\`a della Svizzera italiana; rolf.krause@kaust.edu.sa}}

\maketitle

\begin{abstract}
Training deep neural networks at scale can benefit from domain decomposition, where the network is split into subdomains trained in parallel and coupled by a global trust-region mechanism. Building on the Additively Preconditioned Trust-Region Strategy (APTS), we propose a non-monotone variant with a nonlinear additive Schwarz preconditioner that combines parallel subdomain corrections with global coarse-space directions. A windowed acceptance criterion allows controlled objective increases, avoiding needless rejection of effective coarse steps. The resulting non-monotone APTS (NAPTS) preserves accuracy while reducing CPU time by 30\% and cutting rejected steps to one third of those in APTS.
\end{abstract}

\section{Introduction}
Given samples \(\{(x^{(i)},y^{(i)})\}_{i=1}^M\) with \(x^{(i)}\in\mathbb{R}^q\) and \(y^{(i)}\in\mathbb{R}^p\), we train a neural network (NN) by solving for \(\theta\in\mathbb{R}^n\)
\begin{equation}\label{def:min-problem}
\min_{\theta\in\mathbb{R}^n} f(\theta)
\;=\;
\frac{1}{M}\sum_{i=1}^M \mathcal{L}\left(\mathcal{N}(\theta;x^{(i)}),\,y^{(i)}\right),
\end{equation}
where \(\mathcal{L}:\mathbb{R}^p\times\mathbb{R}^p\to\mathbb{R}\) is a per-sample loss and \(\mathcal{N}(\theta;\cdot)\) maps inputs to predictions. In modern applications, \(f\) is large-scale, nonconvex, and typically evaluated under sampling noise.

First-order methods such as stochastic gradient descent (SGD) and Adam remain the default \cite{Goodfellow-Bengio-Courville-2016,Kingma-Ba-2014}, with performance largely governed by step-size and momentum schedules. This leads to many hyper-parameters that must be tuned and does not exploit the pipelined execution of the model across multiple devices, where a NN is split into sequential blocks placed on different GPUs and activations are forwarded from one device to the next. These issues motivate the study of algorithms that explicitly exploit such multi-device pipelines, and rely on globalization strategies which adapt parameters (e.g., step sizes) automatically.

Domain decomposition (DD) enables model-parallelism by partitioning parameters into subdomains and updating restricted subproblems in parallel with additive recombination\cite{Chan-Zou-1994,Toselli-Widlund-2004, Erhel-Gander-Halpern-Pichot-Sassi-Widlund-2014}.
This mirrors classical DD for PDEs \cite{Toselli-Widlund-2004} and aligns with modern model-parallel pipelines \cite{BenNun-Hoefler-2019,Nichols-Singh-Lin-Bhatele-2021}. 
For globalization, trust-region (TR) strategies provide robust control for nonconvex problems \cite{Nocedal-Wright-1999,connGouldToint2000} and have been adapted to multi-level and mini-batch regimes \cite{Curtis-Scheinberg-Shi-2019,Kopanicakova-Krause-2022}. Within DD, additive updates with a TR safeguard underlie the Additively Preconditioned Trust-Region Strategy (APTS)\cite{Gross-Thesis}.

In this paper, we extend the APTS framework \cite{CruzAlegria-Capriqi-Likaj-Trotti-Krause-2025-arXiv, Gross-Thesis,CruzAlegria-et-al-2025,Trotti-et-al-2025} by incorporating the principles of non-monotone TR (NTR) methods, resulting in the Non-monotone Additively Preconditioned Trust-Region method (NAPTS). NAPTS compares trial steps against a sliding reference value, which relaxes strict monotonic decrease, accommodates noisy mini-batch objectives, and enables more persistent coarse-space directions.

\section{Foundations of NAPTS} 

\label{sec:theory}
This section describes the subdomains structure, the local solves, and the globalization strategy.
\subsection{Additive domain decomposition update} 
\label{subsec:additive-minimal}
Let $\{C_d\}_{d=1}^N$ be a partition of $\{1,\dots,n\}$ with restriction and prolongation operators $R_d^{\phantom{\mathrm{T}}}:\mathbb{R}^n\to\mathbb{R}^{n_d}$ and
$R_d^\mathrm{T}:\mathbb{R}^{n_d}\to\mathbb{R}^n$ satisfying
$R_d^{\phantom{\mathrm{T}}} R_d^\mathrm{T} = I_{n_d}^{\phantom{\mathrm{T}}}$, $R_d^{\phantom{\mathrm{T}}} R_{d'}^\mathrm{T}=0$ for $d\neq d'$, and
$\sum_{d=1}^N R_d^\mathrm{T} R_d^{\phantom{\mathrm{T}}} = I_n^{\phantom{\mathrm{T}}}$. This is the non overlapping parameter space analogue of an additive Schwarz decomposition, where the ``domain'' is the NN parameter space.

\begin{figure}[h!]
    \centering
    \includegraphics[width=0.7\linewidth]{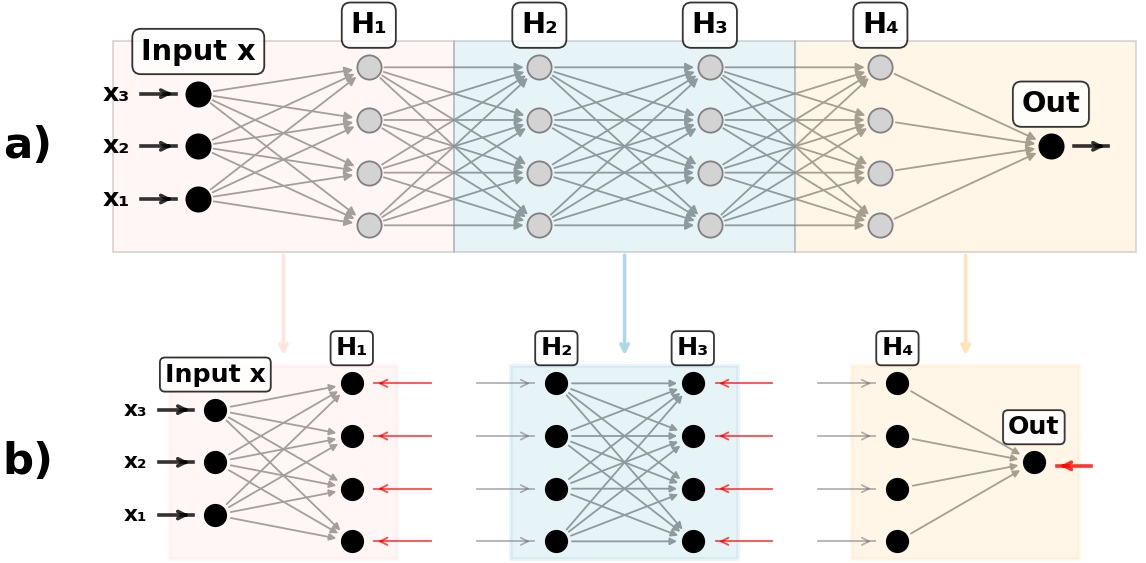}
    \caption{A NN decomposition example. (a) Hardware-based partition into three blocks. (b) Black/red arrows denote cached boundary activations/downstream gradients.}
    \label{fig:nnpg}
\end{figure}

To make this concrete, Figure~\ref{fig:nnpg}a shows an example of such a decomposition with $N=3$ for a fully connected NN, where each color identifies one subdomain of the parameter partitions $C_d$. 

To let the subdomains operate independently, we first expose them to global information. Therefore, given a global iterate $\theta^k$ at iteration $k$ of NAPTS, we compute the objective value $f{\left(\theta^k\right)}$ and its gradient $\nabla f\left(\theta^k\right)$ by running a single forward and backward propagation through the full NN, during which we cache the activation values and downstream gradients  (black/red arrows in Figure \ref{fig:nnpg}b). To illustrate this idea for the simple case of the NN in Figure~\ref{fig:nnpg}, let us define the subdomains $D_d$ as:
\[
\begin{array}{r l}
\text{Subdomains:}\qquad &
\begin{aligned}
D_1(x;\theta_1) &= \sigma_1(W_{H_1} x + b_{H_1}) ,\\
D_2(x;\theta_2) &= \sigma_2\bigl(W_{H_2}\,\sigma_3(W_{H_3}x+b_{H_3}) + b_{H_2}\bigr),\\ 
D_3(x;\theta_3) &= \sigma_{\text{out}}\left(W_{\text{out}}\,\sigma_4 \left(W_{H_4}x+b_{H_4}\right) + b_{\text{out}}\right),
\end{aligned}
\\[0.6em]
\text{Full NN:}\qquad &
\mathcal{N}(x;\theta) = D_3(D_2(D_1(x))),
\end{array}
\]
where $\sigma_i$ are the activation functions and $\theta$ contains all the network parameters $\theta_d= R_d \theta \    \forall    d $, which split into weights $W$ and biases $b$.
By the chain rule and linearity, the gradient $\nabla f(\theta)$ can be
factored into a downstream gradient and a local derivative, as
summarized in Table~\ref{tab:block-gradients} where $\hat y$ denotes the
NN output on the full dataset or current minibatch. 

\begin{table}[h]
\centering
\renewcommand{\arraystretch}{2}
\begin{tabular}{|c|c|c|} 
\hline
\textbf{Parameter $\theta_d$} &
\textbf{Downstream gradient $G_d$} &
\textbf{Gradient $\dfrac{\partial f}{\partial \theta_d}$} \\[0.8ex] 

\hline
$\theta_3$ &
$G_3 = \dfrac{\partial f}{\partial \hat y}
      = \dfrac{\partial f}{\partial D_3}$ &
$\dfrac{\partial f}{\partial \theta_3}
  = G_3 \dfrac{\partial D_3}{\partial \theta_3}$ \\[2mm]
$\theta_2$ &
$G_2 = \dfrac{\partial f}{\partial \hat y}
      \dfrac{\partial \hat y}{\partial D_2}
    = \dfrac{\partial f}{\partial D_2}$ &
$\dfrac{\partial f}{\partial \theta_2}
  = G_2 \dfrac{\partial D_2}{\partial \theta_2}$ \\[2mm]
$\theta_1$ &
$G_1 = \dfrac{\partial f}{\partial \hat y}
      \dfrac{\partial \hat y}{\partial D_2}
      \dfrac{\partial D_2}{\partial D_1}
    = \dfrac{\partial f}{\partial D_1}$ &
$\dfrac{\partial f}{\partial \theta_1}
  = G_1 \dfrac{\partial D_1}{\partial \theta_1}$ \\[1mm]
\hline
\end{tabular}
\caption{Stored downstream and blockwise gradients for each $\theta_d$.}
\label{tab:block-gradients}
\end{table}

Thus, during local iterations, the cached forward activation from the previous subdomain provides the input needed to evaluate the local derivative 
{$\partial D_d / \partial \theta_d$}, while the stored backward gradient from the next subdomain, $G_{d+1}=\partial f/\partial D_{d+1}$, allows us to reconstruct $G_d$. Starting from $\theta_d^{k,0} = R_d \theta^k$,
during the $\ell$ inner iterations $t=1,...,\ell$, we keep the stored downstream gradient $G_d^k$ fixed,
while recomputing the local derivative
$\partial D_d / \partial \theta_d$ using cached boundary activations. This yields approximate block–gradients of the form
\begin{equation}\label{eq:g_tilde}
  \tilde g_d^{k,t}
  \;:=\;
  G_d^k \,
  \frac{\partial D_d}{\partial \theta_d}\bigl(x_d^k;\,\theta_d^{k,t}\bigr)
  \;\approx\; \frac{\partial f}{\partial \theta_d}\left(\theta^{k,t}\right),
  \qquad t = 1,\dots,\ell,
\end{equation}
where $x_d^k$ denotes the cached input of subdomain $d$ at the iteration $k$.

For the local subdomain updates, we use an Adam-type method where the gradients are computed through equation \eqref{eq:g_tilde} and Adam step is truncated to satisfy the per-step bound $\Delta_d^{k} := \Delta^k / \ell$, i.e.,
\begin{equation}\label{eq:ClipAdamstep}
s_d^{k,t} =\min \left\{ \left\|\tilde s_d^{k,t}\right\|_{\infty},\Delta_d^{k} \right\} \dfrac{\tilde{s}_d^{k,t}}{\left\|\tilde s_d^{k,t}\right\|_{\infty}}  \quad \text{and update \quad \(
  \theta_d^{k,t} \;=\; \theta_d^{k,t-1} +  s_d^{k,t}.
\)
} 
\end{equation}
Therefore, after $\ell$ local iterations, the subdomain $d$ yields
\(s_d^k=\sum_{t=1}^{\ell}  s_d^{k,t},
  \) such that \( \|s_d^k\|_\infty \le \Delta_d^k.
\)
The global trial step is the additive lift of the $N$ local updates $\{s_d^k\}_{d=1}^N$,
\begin{equation}\label{eq:additive-global-step}
s^k \;=\; \sum_{d=1}^N R_d^\mathrm{T} s_d^k.
\end{equation}

Finally, we note that communications are not needed during the subdomain training and that this construction is not tied to the NN structure in Figure~\ref{fig:nnpg} and generalizes to more complex architectures.

\subsection{Non-monotone trust-region as globalization strategy}\label{sec:NTR} 
In this section, we follow the NTR framework~\cite[Alg.~10.1.1]{connGouldToint2000}. 
Given the proposed additive search direction $s^k$ in~\eqref{eq:additive-global-step}, the NTR mechanism acts, in NAPTS, as a globalization strategy deciding whether to accept the global step $s^k$. In contrast to previous APTS variants~\cite{CruzAlegria-Capriqi-Likaj-Trotti-Krause-2025-arXiv,CruzAlegria-et-al-2025,Trotti-et-al-2025} with exact gradient and based on monotone TR schemes, we deliberately switch to a non-monotone variant to increase the likelihood of accepting coarse space steps, and thereby avoid expensive recomputation of the search direction. This is particularly appropriate in our setting, where the coarse-space updates are inexact and potentially batch-based, and are therefore inherently non-monotone. Preliminary experiments in~\cite{Trotti-et-al-2025} suggest that aggressively accepting coarse-space directions can accelerate loss decrease even at the price of temporary objective increases, but this comes at the cost of an indiscriminate acceptance policy. In NAPTS, we therefore introduce a more robust NTR mechanism that selectively accepts non-monotone coarse-space steps in a controlled way.

We introduce a local model $m^k$ that approximates $f$ in the global trust-region
\[
  \mathcal B_G^k \;:=\; \bigl\{\,s\in\mathbb{R}^n : \|s\|_\infty \le \Delta^k \bigr\},
  \qquad \Delta^k>0.
\]
Classically, TR methods employ a quadratic model and define the trial step as an approximate solution of the corresponding trust-region subproblem over $\mathcal B_G^k$. However, in our large-scale setting, forming and manipulating exact second-order information is prohibitive, and using approximate second-order terms introduces additional complications. For simplicity, we therefore focus on a first-order model
\begin{equation}\label{eq:TR-model-firstorder}
  m^k(s)
  \;=\;
  f\left(\theta^k\right)
  \;+\;
  \nabla f\left(\theta^k\right)^\mathrm{T} s ,
\end{equation}
and use step $s^k$ in equation~\eqref{eq:additive-global-step} as trial step that satisfies $\left\|s^k\right\|_\infty\le \Delta^k$.

Let $\nu\in\mathbb{N}$ be the memory parameter and let $\mathcal W_\nu^k$ be the set that keeps track of the indices of the most recent $\nu$ successful iterations, i.e.
\begin{equation}\label{eq:Wset}
  \mathcal W_\nu^k 
  \;:=\;
  \bigl\{\, j \in \{ \max\{0,k-\nu\},\dots,k \} 
  \;\big|\; \text{iteration } j \text{ is successful} \,\bigr\}.
\end{equation}
We then define the reference-index map $r:\mathbb{N} \to \mathbb{N}$ as
\begin{equation}\label{eq:def-r-argmax}
  r(k)
  \;:=\;
  \argmax_{j \in \mathcal{W}^k_\nu}
  f\left(\theta^j\right),
\end{equation}
so that  $r(k)\in \mathcal{W}^k_\nu$. The corresponding history term 
$\sigma^k_h \ge 0$ is defined as the cumulative predicted expected decrease among successful 
iterations between the reference index $r(k)$ and the current one,
\begin{equation}\label{eq:sigma-def}
  \sigma^k_h
  \;:=\;
  \sum_{\substack{i = r(k) \\ i \in \mathcal W_\nu^k}}^{k-1}
  \bigl[m^i(0) - m^i( s^i)\bigr].
\end{equation}

Once $s^k$ is computed, the quality of the model prediction is assessed via two agreement ratios
\begin{equation}\label{eq:rho}
  \rho^k_c
  \;=\;
  \frac{f\left(\theta^k\right)-f\left(\theta^k+s^k\right)}
       {m^k\left(0\right)-m^k\left(s^k\right)},
  \qquad
  \rho^k_h
  \;=\;
  \frac{f\left(\theta^{r(k)}\right)-f\left(\theta^k+s^k\right)}
       {\sigma_h^k+m^k\left(0\right)-m^k\left(s^k\right)},
\end{equation}
where the subscripts $c$ and $h$ refer to \textit{current} and \textit{historical} quantities,
respectively.
Specifically, we set \(  \rho^k \;:=\; \max\left\{\rho^k_c,\rho^k_h\right\}  \)
and use this combined ratio in place of the classical TR ratio $\rho^k_c$ to decide step acceptance and to update the global radius $\Delta^k$.

By construction, $f\left(\theta^{r(k)}\right) \ge f\left(\theta^k\right)$, so $\rho_h^k$ may be
larger than $\rho_c^k$ and accept the step even in the case where $f\left(\theta^k+s^k\right) > f\left(\theta^k\right)$. Hence, the sequence
$\left\{f\left(\theta^k\right)\right\}_k$ is not required to be monotonically decreasing, while the
reference values $\left\{f\left(\theta^{r(k)}\right)\right\}_k$ and the history term $\sigma_h^k$ still enforce a long-term descent behaviour. The additional overhead with respect to a standard monotone TR scheme is negligible: maintaining the window
$\mathcal W_\nu^k$, the index $r(k)$, and the scalar $\sigma_h^k$ requires only 
simple updates per iteration, yet this non-monotone mechanism can substantially
increase the number of accepted steps coming from coarse spaces, thereby
reducing the need for expensive global corrections of the search direction.

\subsection{The NAPTS Method}
\label{sec:apts}
The NAPTS method, summarized in Algorithm~\ref{alg:apts}, combines parallel subdomain updates with two NTR stages. Each iteration naturally decomposes into two phases: parallel subdomain computations, preconditioning, and NTR steps.
In the first phase, lines~\ref{line:local-for}–\ref{line:additive}, each subdomain $D_d$ performs in parallel $\ell$ local constrained Adam(CAdam) steps. This produces local steps $s_d^k$ which are lifted and combined to form the global proposal $s^k$ as in \eqref{eq:additive-global-step}. 

The second phase (lines~\ref{line:ntr-history}–\ref{line:ntr-step}) applies the NTR update strategy to the obtained proposal $s^k$. Using the reference quantities $\mathcal W_\nu^k$, $r(k)$, and $\sigma_h^k$, we compute the current and historical agreement ratios $\rho_c^k$ and
$\rho_h^k$ and combine them into $\rho^k = \max\{\rho_c^k,\rho_h^k\}$
(cf.~\eqref{eq:rho}). $\rho^k$ is then used to decide whether to take the step $s^k$ or a corrected step $c^k$, and to update the global radius $\Delta^k$. The correction $c^k$ 
in line~\ref{line:c_k} is a convex combination of a steepest descent step of 
length $\Delta^k$ and the global proposal $s^k$, for some 
$\alpha_k, \beta_k \in [0,1]$, namely,
\begin{equation}\label{eq:c_k}  
  c^k \;=\;\beta_k\left[ (1-\alpha_k)\Bigl(-\Delta^k \frac{\nabla f(\theta^k)}{\|\nabla f(\theta^k)\|}\Bigr)
            \;+\; \alpha_k s^k\right].
\end{equation}

Coefficients $(\alpha,\beta)\in\left\lbrace (0.8,\frac{1}{2}), (0.6,\frac{1}{4}), (0.4,\frac{1}{8}), (0.2,\frac{1}{16}), (0,\frac{1}{32})\right\rbrace$ are tested through a for loop. This construction avoids wasting the costly subspace updates and gradient computations when the proposed step $s^k$ would otherwise be rejected. 

In the third phase at line~\ref{line:ntr-pure}, a classical NTR step is performed, as described in Section~\ref{sec:NTR}, starting from $(\theta^{k+\frac{1}{2}},\Delta^{k+\frac{1}{2}})$ and using a search direction based on a first-order model. We note that setting $\nu\!=\!1$ recovers the APTS method studied in~\cite{Trotti-et-al-2025} with approximated subdomains.

\begin{algorithm}[t]
\footnotesize
\caption{\footnotesize NAPTS: Non-monotone Additively Preconditioned Trust-Region Strategy}
\label{alg:apts}
\begin{algorithmic}[1]
\Require {$f:\R^n \to \R$, initial iterate $\theta^0 \in \R^n$, restriction operators $\{R_d\}_{d=1}^N$ }
\Ensure Approximate minimizer $\theta^\ast$ of $f$
\Constants {inner iterations $\ell$, memory $\nu$, NTR parameters $(\eta_1,\eta_2,\gamma_{\mathrm{dec}},\gamma_{\mathrm{inc}},\Delta^0)$}\Comment{See \cite{connGouldToint2000}}
\State Set $k \gets 0$ and initialize NTR history $(\sigma_h^0,r(0),\mathcal W_\nu^0)$ \Comment{cf. Section~\ref{sec:NTR}}
\While{not converged}
  \State Compute $f(\theta^k)$, $\nabla f(\theta^k)$, block gradients $g_d^k \gets R_d \nabla f(\theta^k)$
 \ForAll{$d = 1,\dots,N$ \textbf{ in parallel}} \label{line:local-for} \Comment{\textit{Phase 1:} parallel subdomains}
           \State \hspace{1.5em}$s_d^k \gets \operatorname{CAdam}\bigl(g_d^k, \ell, \Delta_d^k\bigr)$ \Comment{local step, cf.~\eqref{eq:ClipAdamstep}}
         \EndFor
  \State \label{line:additive} Form the global proposal $s^k = \sum_{d=1}^N R_d^\mathrm{T} s_d^k$ 
 
\State \label{line:jump1}Compute $\rho_c^k$, $\rho_h^k$, $\rho^k=\max\{\rho_c^k,\rho_h^k\}$  \Comment{\textit{Phase 2:} NTR, cf.~\eqref{eq:rho}}\label{line:ntr-history}
  \State Step update\label{line:c_k}\Comment{cf.~\eqref{eq:c_k}}\vspace{-4mm}
  \[
    \theta^{k+\frac{1}{2}} :=
    \begin{cases}
      \theta^{k} + s^k, & \text{if } \rho^k > \eta_1,\\
      \theta^{k} +c^k,       & \text{otherwise,}
    \end{cases}
    \qquad
    \Delta^{k+\frac{1}{2}} :=
    \begin{cases}
      \gamma_{\mathrm{inc}}\Delta^k, & \text{if } \rho^k \ge \eta_2,\\
      \Delta^k,                      & \text{if } \rho^k \in [\eta_1,\eta_2),\\
      \gamma_{\mathrm{dec}}\Delta^k, & \text{if } \rho^k < \eta_1.
    \end{cases}
  \]\vspace{-4mm}
   \State Update\label{line:ntr-step} $\mathcal{W}_\nu^k$, $r(k)$, $\sigma_h^k$, \Comment{ cf.~\eqref{eq:Wset}–\eqref{eq:sigma-def}}
   \State \label{line:ntr-pure}$\theta^{k+1},\Delta^{k+1}\gets\operatorname{NTR}(\theta^{k+\frac{1}{2}},\Delta^{k+\frac{1}{2}},\nabla f(\theta^{k+\frac{1}{2}}))$\Comment{\textit{Phase 3:} Smoothing NTR iteration}
  \State $k\gets k+1$
\EndWhile
\State \Return $\theta^\ast := \theta^k$
\end{algorithmic}
\end{algorithm}

\section{Numerical examples}\label{Sec:NumericalEX}
We consider image classification on the \texttt{CIFAR-10} dataset. 
For this task, we use a convolutional NN (CNN\footnote{A CNN is a feedforward architecture that applies learned convolutional filters; see \cite[Ch.~9]{Goodfellow-Bengio-Courville-2016}.}) with four convolutional blocks and two fully-connected layers (1.2M parameters) trained in a multi-GPU distributed setting on a compute node with four Nvidia A100 GPUs using data mini-batches of size 1{,}000. The model decomposition is layer-based, with each subdomain corresponding to a contiguous block of layers, consistent with the schematic illustrated in Figure~\ref{fig:nnpg}.

\begin{figure}
    \centering
    \includegraphics[width=1\linewidth]{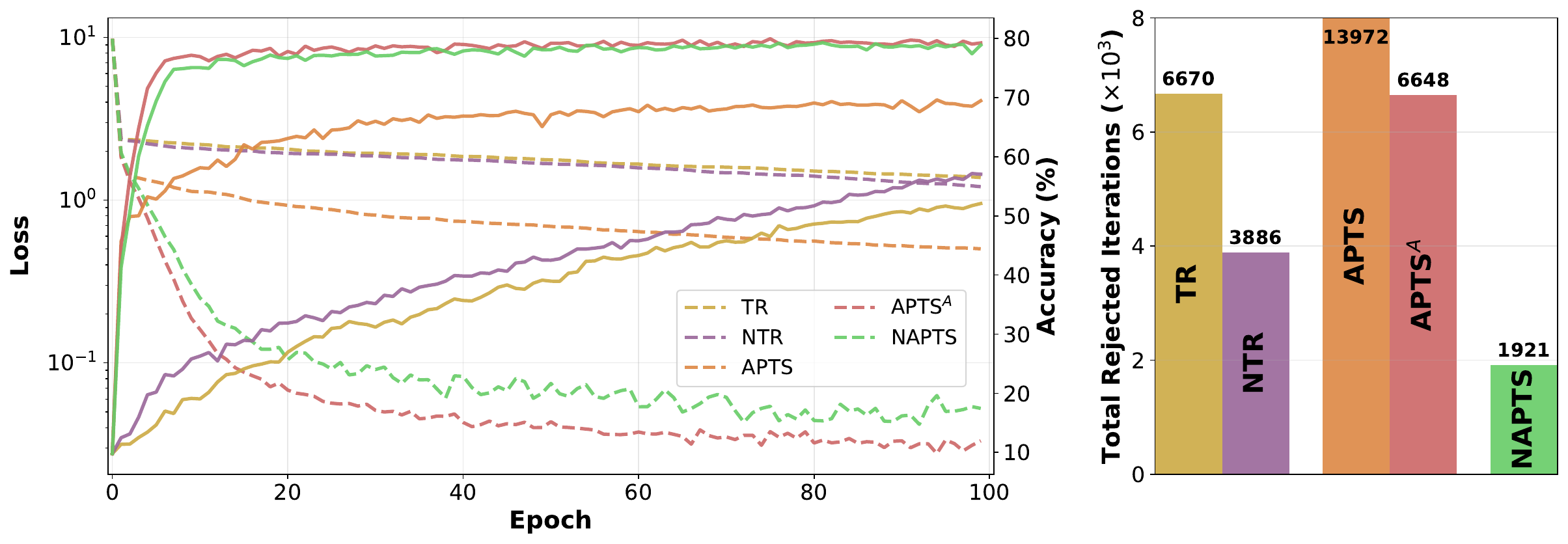}
    \caption{Loss and accuracy (left) and rejected steps per batch (right). NTR and NAPTS use $\nu=100$, and the superscript $A$ denotes always accepting the global step in~\eqref{eq:additive-global-step}.}
    \label{fig:results}
\end{figure}

A rejection is counted whenever a search direction fails a TR acceptance test, either at the globalization level or within the internal correction loop used to construct $c^k$.

Left panel of Figure~\ref{fig:results} reports the training loss and validation accuracy, and the right panel the rejection counts per batch, for TR (53.65s/epoch), NTR (47.26s/epoch), APTS (169.14s/epoch), APTS$^{A}$ (152.39s/epoch), and the proposed method NAPTS (115.88s/epoch).
Training loss is an optimization diagnostic and does not, by itself, imply better predictive performance. For this reason, we assess performance primarily through the validation accuracy, i.e., the fraction of correctly classified samples on the held-out validation set. 

The baseline TR and NTR schemes struggle to reach high accuracy. Although they require roughly half the backpropagations per epoch relative to the APTS-based methods, they still fall short of comparable accuracy even after twice as many epochs. Nevertheless, NTR improves over TR, with about half the rejected directions and better loss and accuracy. This gap motivates a non-monotone globalization within APTS, which is well suited to the nonconvex and noisy training landscape.

In contrast, the preconditioned variants APTS and NAPTS achieve lower training losses and substantially higher validation accuracies, suggesting improved optimization without loss of generalization. As shown in \cite{Trotti-et-al-2025}, the aggressive-acceptance variant APTS$^A$ is also effective, yielding about a $10\%$ reduction in CPU time relative to APTS by always accepting the global step in~\eqref{eq:additive-global-step}. Moreover, NAPTS provides a clear improvement over APTS$^A$ by retaining a controlled loss-based acceptance mechanism that enhances robustness, while delivering a further $\approx 20\%$ reduction in CPU time compared to APTS$^A$ (about $30\%$ relative to APTS). This gain follows directly from the marked decrease in rejected steps visible in the right panel, indicating that the preconditioned step is almost always accepted. 

In our setup, Adam/SGD requires $\approx 40$s/epoch (about $3\times$ faster), but typically needs learning-rate and schedule tuning across multiple runs. In contrast, APTS/NAPTS uses trust regions to adapt step sizes automatically. Overall, NAPTS improves acceptance of coarse-space proposals, reducing rejections and CPU time.

\section*{Acknowledgements}
The research was funded by the Swiss National Science Foundation (SNSF) under projects No. 224943 and No. 197041 (ML$^2$).
\bibliographystyle{spmpsci}
\bibliography{Biblio}      

@inproceedings{CruzAlegria-et-al-2025,
author = {{Cruz Alegr{\'i}a}, S. and Trotti, K. and Kopani{\v{c}}{\'a}kov{\'a}, A. and Krause, R.},
  title     = {Data-Parallel Neural Network Training via Nonlinearly Preconditioned Trust-Region Method},
booktitle = {ENUMATH 2023},
  series    = {Lect. Notes Comput. Sci. Eng.},
  volume    = {153},
  pages     = {34--43},
  year      = {2025},
  address   = {Berlin, Germany},
  publisher = {Springer}
}

@incollection{Trotti-et-al-2025,
  author = {Trotti, K. and {Cruz Alegr{\'i}a}, S. and Krause, R. and Kopani{\v{c}}{\'a}kov{\'a}, A.},
  title     = {Parallel Trust‐Region Approaches in Neural Network Training},
  booktitle = {Proceedings of the MATH+ Thematic Einstein Semester 2023:\\Mathematical Optimization for Machine Learning},
  publisher = {De Gruyter},
  address   = {Berlin},
  pages     = {107--120},
  year      = {2025},
  isbn      = {9783111376776},
}

@book{connGouldToint2000,
	title        = {{Trust region methods}},
	author = {Conn, A. R. and Gould, N. I. and Toint, P. L.},
	year         = {2000},
	publisher    = {Society for Industrial and Applied Mathematics}
}

@inproceedings{Kingma-Ba-2014,
  author       = {Kingma, D. P. and Ba, J.},
  title        = {Adam: A Method for Stochastic Optimization},
  booktitle    = {3rd International Conference on Learning Representations (ICLR 2015), San Diego, CA, USA, May 7--9, 2015},
  year         = {2015},
  eprint       = {1412.6980},
  archivePrefix= {arXiv},
  primaryClass = {cs.LG}
}

@book{Nocedal-Wright-1999,
  author    = {Nocedal, J. and Wright, S.J.},
  title     = {{Numerical Optimization}},
  publisher = {Springer},
  address   = {New York, NY},
  year      = {1999},
}

@phdthesis{Gross-Thesis,
  author  = {Gro{\ss}, C.},
  title   = {{A Unifying Theory for Nonlinear Additively and Multiplicatively Preconditioned Globalization Strategies: Convergence Results and Examples From the Field of Nonlinear Elastostatics and Elastodynamics}},
  school  = {Bonn International Graduate School, University of Bonn},
  address = {Bonn, Germany},
  year    = {2009},
}

@article{Curtis-Scheinberg-Shi-2019,
  author = {Curtis, F. E. and Scheinberg, K. and Shi, R.},
  title   = {{A Stochastic Trust-Region Algorithm Based on Careful Step Normalization}},
  journal = {INFORMS J. Optim.},
  volume  = {1},
  pages   = {200--220},
  year    = {2019},
}

@article{Kopanicakova-Krause-2022,
  author = {Kopani{\v{c}}{\'a}kov{\'a}, A. and Krause, R.},
  title   = {{Globally Convergent Multilevel Training of Deep Residual Networks}},
  journal = {SIAM J. Sci. Comput.},
  volume  = {45},
  number  = {3},
  pages   = {S254--S280},
  year    = {2023},
}

@book{Toselli-Widlund-2004,
  author    = {Toselli, A. and Widlund, O.},
  title     = {{Domain Decomposition Methods: Algorithms and Theory}},
  series    = {Springer Ser. Comput. Math.},
  volume    = {34},
  publisher = {Springer},
  address   = {Berlin, Germany},
  year      = {2004},
}

@misc{Nichols-Singh-Lin-Bhatele-2021,
  author = {Nichols, D. and Singh, S. and Lin, S. H. and Bhatele, A.},
  title        = {{A Survey and Empirical Evaluation of Parallel Deep Learning Frameworks}},
  howpublished = {arXiv preprint arXiv:2111.04949},
  year         = {2021},
}

@article{BenNun-Hoefler-2019,
  author = {Ben-Nun, T. and Hoefler, T.},
  title   = {{Demystifying Parallel and Distributed Deep Learning: An In-Depth Concurrency Analysis}},
  journal = {{ACM} Comput. Surv.},
  volume  = {52},
  number  = {1},
  pages   = {1--43},
  year    = {2019},
}

@article{Chan-Zou-1994,
  author = {Chan, T. F. and Zou, J.},
  title   = {{Additive Schwarz Domain Decomposition Methods for Elliptic Problems on Unstructured Meshes}},
  journal = {Numer. {Algorithms}},
  volume  = {8},
  number  = {2},
  pages   = {329--346},
  year    = {1994},
}

@book{Erhel-Gander-Halpern-Pichot-Sassi-Widlund-2014,
  author = {Erhel, J. and Gander, M. J. and Halpern, L. and Pichot, G. and Sassi, T. and Widlund, O.},
  title     = {{Domain Decomposition Methods in Science and Engineering XXI}},
  publisher = {Springer},
  address   = {Switzerland},
  year      = {2014},
}

@book{Goodfellow-Bengio-Courville-2016,
  author = {Goodfellow, I. and Bengio, Y. and Courville, A.},
  title         = {{Deep Learning}},
  publisher     = {MIT Press},
  year          = {2016}
}

@misc{CruzAlegria-Capriqi-Likaj-Trotti-Krause-2025-arXiv,
  author        = {Cruz Alegr{\'i}a, S. and {\c{C}}apriqi, B. and Likaj, S. and Trotti, K. and Krause, R.},
  title         = {{An Additively Preconditioned Trust-Region Strategy for Machine Learning}},
  howpublished = {arXiv preprint arXiv:2512.14286},
  year         = {2025},
}

\end{document}